\documentclass[amsbook,12pt]{article}
\usepackage{epsfig, amsfonts,amssymb, amsmath}
\usepackage[usenames]{color}
\usepackage{picinpar}
\usepackage{float, graphicx}
\newcommand{\nl}{\mbox{}\\}

\begin{document}
%
%
%
\mbox{} \vspace{-2.000cm} \\
\begin{center}
\mbox{\Large \bf %
On the global solvability of porous media equations} \\
\mbox{} \vspace{-0.250cm} \\
\mbox{\Large \bf %
with general (spatially dependent) advection terms} \\
\nl
\mbox{} \vspace{-0.200cm} \\
{\large \sc N.\;M.\;L.\;Diehl,}$\mbox{}^{\!\:\!1}$
{\large \sc L.\;Fabris}$\mbox{}^{\;\!2}$
{\large \sc and P.\;R.\;Zingano}$\mbox{}^{\;\!3}$ \\
\mbox{} \vspace{-0.175cm} \\
$\mbox{}^{1}${\small
Instituto Federal de Educa\c c\~ao, Ci\^encia e Tecnologia} \\
\mbox{} \vspace{-0.685cm} \\
{\small
Canoas, RS 92412, Brazil} \\
\mbox{} \vspace{-0.390cm} \\
$\mbox{}^{2}${\small
Coordenadoria Acad\^emica} \\
\mbox{} \vspace{-0.685cm} \\
{\small
Universidade Federal de Santa Maria\,-\;\!\:\!Cachoeira do Sul} \\
\mbox{} \vspace{-0.685cm} \\
{\small
Cachoeira do Sul, RS 96501, Brazil} \\
\mbox{} \vspace{-0.390cm} \\
$\mbox{}^{3}${\small Departamento de Matem\'atica Pura e Aplicada} \\
\mbox{} \vspace{-0.670cm} \\
{\small Universidade Federal do Rio Grande do Sul} \\
\mbox{} \vspace{-0.670cm} \\
{\small Porto Alegre, RS 91509, Brazil} \\
\mbox{} \vspace{-0.390cm} \\
\nl
%
%
%
%
%
\mbox{} \vspace{-0.400cm} \\
{\bf Abstract} \\
\mbox{} \vspace{-0.525cm} \\
\begin{minipage}[t]{12.250cm}
{\small
\mbox{} \hspace{+0.150cm}
We show that
advection-diffusion equations
with porous media type diffusion
and integrable
initial data
are globally solvable
under very mild assumptions.
Some generalizations
and related results
are also given. \\
}
\end{minipage}
\end{center}
%
%
%
%
\nl
\mbox{} \vspace{-0.650cm} \\
\mbox{} \hspace{+0.250cm}
{\small 35A01},
{\small 35K15} \\
\mbox{} \vspace{-0.250cm} \\
\mbox{} \hspace{+0.250cm}
\mbox{} \hspace{+0.850cm}
\mbox{} \hspace{+0.100cm}
Cauchy problem,
spatially dependent advection flux,
pointwise estimates \\
\nl
\setcounter{page}{1}
\mbox{} \vspace{-0.250cm} \\
%
%
%
%

%
%

%
{\bf 1. Introduction} \\
\mbox{} \vspace{-0.650cm} \\

In this note,
we describe general results
recently obtained by the authors
concerning
the solvability in the large
of initial value problems
for
degenerate advection-diffusion
equations
of the type \\
\mbox{} \vspace{-0.650cm} \\
\begin{equation}
\tag{1.1$a$}
u_t \,+\: \mbox{div}\;\!\mbox{\boldmath $f$}(x,t,u)
\,+\: \mbox{div}\,\mbox{\boldmath $g$}(t,u)
\:=\;
\mu(t) \;\!\;\! \mbox{div}\,(\;\! |\,u\,|^{\:\!\alpha}
\,\nabla u \;\!),
\end{equation}
\mbox{} \vspace{-0.950cm} \\
\begin{equation}
\tag{1.1$b$}
u(\cdot,0) \,=\,
u_{0} \in L^{1}(\mathbb{R}^{n}) \cap
L^{\infty}(\mathbb{R}^{n})
\end{equation}
\mbox{} \vspace{-0.250cm} \\
and some generalizations
(see Section~2).
Here,
$ \alpha > 0 $ is constant,
$ \mu \in C^{0}(\:\![\;\!0, \infty)\,\!) $
is positive,
and
$ \mbox{\boldmath $f$}\!\:\!= (\:\!f_{\mbox{}_{1}} \!\:\!,
f_{\mbox{}_{2}} \!\;\!,\!\;\!...\:\!, f_{\mbox{}_{\scriptstyle n}}) $,
$ \mbox{\boldmath $g$}\!\:\!= (\;\!g_{\mbox{}_{1}} \!\!\;\!\;\!,
g_{\mbox{}_{2}} \!\;\!,\!\;\!...\:\!, g_{\mbox{}_{\scriptstyle n}}) $,
%
%
are
given continuous
advection flux fields
that are locally Lipschitz in $u$
uniformly in $ \:\!x \in \mathbb{R}^{n} \!\;\!$
and bounded $\:\! t \geq 0 $,
with
%
%
$ \!\;\!\mbox{\boldmath $f$} \!\;\! $
satisfying\:\!:
$ \!\;\!\mbox{\boldmath $f$}(x,t,0) = {\bf 0} $
for all $ \:\!x, \:\!t \;\!$
and \\
\mbox{} \vspace{-0.675cm} \\
\begin{equation}
\tag{1.2}
|\,\mbox{\boldmath $f$}(x,t,\mbox{u}) \,|
\,\leq\,
\mbox{\small $F$}(t) \,|\,\mbox{u}\,|^{\:\!\kappa\,+\,1}
\quad \;\;\,
\forall \; x \in \mathbb{R}^{n} \!,
\; t \geq 0,
\; \mbox{u} \in \mathbb{R}
\end{equation}
\mbox{} \vspace{-0.250cm} \\
for some
$ \mbox{\small $F$} \in C^{0}(\!\;\!\;\![\;\!0, \infty)\,\!) $
and some constant $ \:\!\kappa \geq 0 $,
where
$ \:\!|\!\;\!\;\!\cdot\!\;\!\;\!| $
denotes the absolute value (in case of scalars)
or the Euclidean norm (in case of vectors),
as in~(1.1$a$). \linebreak
By a (bounded) {\em solution\/}
of the problem (1.1) in some time interval
$ [\;\!0, \;\!\mbox{\small $T$}_{\!\ast}) $
is meant \linebreak
any function
$ {\displaystyle
\;\!
u(\cdot,t) \in
C^{0}([\;\!0, \;\!\mbox{\small $T$}_{\!\ast}),
\!\;\!\;\!
L^{1}_{\mbox{\scriptsize loc}}(\mathbb{R}^{n}))
\cap
L^{\infty}_{\mbox{\scriptsize loc}}([\;\!0, \;\!\mbox{\small $T$}_{\!\ast}),
\!\;\!\;\!
L^{1}(\mathbb{R}^{n}) \!\;\!\cap\!\;\! L^{\infty}(\mathbb{R}^{n}))
} $
having
$ {\displaystyle
|\;\!u(\cdot,t)\;\!|^{\:\!\alpha} \!\;\!\;\!u(\cdot,t)
\in
L^{2}_{\mbox{\scriptsize loc}}((\!\;\!\;\!0, \;\!\mbox{\small $T$}_{\!\ast}),
W^{1,\,2}_{\mbox{\scriptsize loc}}(\mathbb{R}^{n}))
} $
which satisfies (1.1$a$)
in distributional sense
(i.e.,
in
$ {\cal D}^{\;\!\prime}(\;\!\mathbb{R}^{n} \!\times\!\!\;\!\;\!
(\!\;\!\;\!0, \;\!\mbox{\small $T$}_{\!\ast})\,\!) $)
and takes the initial value
$ {\displaystyle
\:\!
u(\cdot,0) =\:\! u_0
\,\!
} $,
see e.g.\;\cite{DaskalopoulosKenig2007, %
Vazquez2007, WuZhaoYinLi2001}.
This says,
in particular,
that
$ u(\cdot,t) \rightarrow u_0 $ in
$ L^{1}_{\mbox{\scriptsize loc}}(\mathbb{R}^{n}) $
as $ t \rightarrow 0 $,
and that,
for every
$ \;\!0 < \mbox{\small $T$} \!\:\!< \mbox{\small $T$}_{\!\ast} $
given, \\
\mbox{} \vspace{-0.650cm} \\
\begin{equation}
\tag{1.3$a$}
\|\, u(\cdot,t) \,
\|_{\mbox{}_{\scriptstyle L^{1}(\mathbb{R}^{n})}}
\,\leq\;
\mbox{\small $M$}_{\mbox{}_{\!1}}\!\;\!(\:\!\mbox{\small $T$}),
\quad \;\;\,
\forall \;\,
0 \:\!\leq\:\! t \:\!\leq\:\! \mbox{\small $T$}\!\;\!,
\end{equation}
\mbox{} \vspace{-0.850cm} \\
\begin{equation}
\tag{1.3$b$}
\|\, u(\cdot,t) \,
\|_{\mbox{}_{\scriptstyle L^{\infty}(\mathbb{R}^{n})}}
\leq\,
\mbox{\small $M$}_{\mbox{}_{\!\infty}}\!\;\!(\:\!\mbox{\small $T$}),
\quad \;\;
\forall \;\,
0 \:\!\leq\:\! t \:\!\leq\:\! \mbox{\small $T$} \!\;\!,
\end{equation}
\mbox{} \vspace{-0.165cm} \\
for some bounds
$ {\displaystyle
\mbox{\small $M$}_{\mbox{}_{\!1}}\!\;\!(\:\!\mbox{\small $T$}),
\;\!
\mbox{\small $M$}_{\mbox{}_{\!\infty}}\!\;\!(\:\!\mbox{\small $T$})
} $
depending on $ \mbox{\small $T$} $
(and the solution $u$ considered).
For the {\em local\/} (in time)
existence
of such solutions,
which are typically obtained
by parabolic regularization
or Galerkin approximations,
see e.g.\;\cite{DaskalopoulosKenig2007, 
Diehl2015, Lions1969, Vazquez2007, WuZhaoYinLi2001}.
From the basic theory,
many interesting solution properties are known;
for example,
one has
$ {\displaystyle
\;\!
u(\cdot,t) \in C^{0}([\;\!0, \mbox{\small $T$}_{\!\ast}),
L^{1}(\mathbb{R}^{n}))
} $
and \\
\mbox{} \vspace{-0.600cm} \\
\begin{equation}
\tag{1.4}
\int_{\mbox{}_{\scriptstyle 0}}^{\;\!T}
\!\!\!\;\!
\int_{\mbox{}_{\scriptstyle \!\;\!\mathbb{R}^{n}}}
\!\!\:\!
|\,u(x,t)\,|^{\:\!2\:\!\alpha} \,
|\, \nabla u(x,t) \,|^{\:\!2}
\:
dx\, dt
\:<\,
\infty
\end{equation}
\mbox{} \vspace{-0.100cm} \\
for every
$ {\displaystyle
\:\!
0 < \mbox{\small $T$} \!<
\mbox{\small $T$}_{\!\ast} \!\;\!
} $,
%
%
see e.g.\;\cite{BrazSchutzZingano2013, Diehl2015, %
Fabris2013, Vazquez2007, WuZhaoYinLi2001}.
\!More importantly to us here,
solutions $ u(\cdot,t) $
decrease monotonically
in $ L^{1}(\mathbb{R}^{n}) $,
%
%
so that,
in particular, \\
\mbox{} \vspace{-0.550cm} \\
\begin{equation}
\tag{1.5}
\mbox{} \hspace{+1.000cm}
\|\, u(\cdot,t) \,
\|_{\mbox{}_{\scriptstyle L^{1}(\mathbb{R}^{n})}}
\leq\;
\|\, u_0 \;\!
\|_{\mbox{}_{\scriptstyle L^{1}(\mathbb{R}^{n})}}
\!\:\!,
\qquad \;\,
\forall \;\;
0 < t < \mbox{\small $T$}_{\!\ast}.
\end{equation}
\mbox{} \vspace{-0.175cm} \\
For all that is presently known,
however,
little has been obtained
regarding the solvability
for large $t$
in the general framework
(1.1), (1.2) above,
except in very special situations.
Thus, for example, when
the flux
$\mbox{\boldmath $f$}(x,t,\mbox{u}) $
does {\em not\/} depend explicitly on $x$,
or, when it does,
if it behaves so as to satisfy
special conditions like \\
\mbox{} \vspace{-0.650cm} \\
\begin{equation}
\mbox{} \hspace{+1.000cm}
\tag{1.6}
\sum_{i\,=\,1}^{n}
\;\!
\mbox{u} \;
\frac{\partial \:\!f_{\scriptstyle i}}
{\partial x_{\scriptstyle i}}
(x,t,\mbox{u})
\:\geq\;0,
\quad \;\;\,
\forall \;\!\;\!\;\!
x \in \mathbb{R}^{n} \!\:\!,
\, t \geq 0,
\, \mbox{u} \in \mathbb{R},
\end{equation}
\mbox{} \vspace{-0.050cm} \\
then
solutions are known to be
globally defined,
with
$ {\displaystyle
\;\!
\|\, u(\cdot,t) \,\|_{L^{q}(\mathbb{R}^{n})}
\!\;\!
} $
monotonically decreasing
for every
$ 1 \leq q \leq \infty $,
see e.g.\;\cite{%
DelPinoDolbeault2003, Diehl2015, %
DiehlFabrisZiebell2018, Fabris2013, %
Porzio2009, Vazquez2007}.
In the absence of (1.6),
however,
things get much more complicated
to analyze.
To see why,
let us consider
by way of illustration
the simple example below.
\!Taking
$ {\displaystyle
\mbox{\boldmath $f$}(x,t,\mbox{u})
=
\mbox{\boldmath $b$}(x,t) \,
|\;\!\mbox{u}\;\!|^{\kappa} \;\! \mbox{u}
\;\!
} $
for some $ \:\!\kappa > 0 $,
and (say)
$ {\displaystyle
\;\!
\mbox{\boldmath $g$}
= {\bf 0}
} $,
$ \mu(t) = 1 $,
the equation (1.1$a$)
becomes \\
\mbox{} \vspace{-0.650cm} \\
\begin{equation}
\tag{1.7}
\notag
u_t \;\!+\,
\;\!
\mbox{\boldmath $b$}(x,t) \cdot
\nabla \:\!(\:\!|\;\!u\;\!|^{\:\!\kappa} \:\!u \;\!)
\:=\;
\mbox{div}\,(\,|\,u\,|^{\:\!\alpha}
\, \nabla u \;\!)
\,+\,
\beta(x,t) \, |\;\!u\;\!|^{\:\!\kappa} \:\!u
\end{equation}
\mbox{} \vspace{-0.250cm} \\
with
$
\:\! \beta(x,t) \:\!=\:\!
- \sum_{\,i\,=\,1}^{\,n}
\!\;\!
\partial \;\!b_{\scriptstyle i}/
\partial \,\!x_{\scriptstyle i}
$.
In regions where
$ \!\;\!\;\!\beta(x,t) \!\;\!> 0 $,
it is clear that
$ \!\;\!\;\!|\!\;\!\;\!u(x,t)\!\;\!\;\!| \!\;\!\;\!$
tends to grow,
particularly if
$ \!\;\!\;\!\beta(x,t) \gg 1 $,
potentially leading to finite-time blow-up. \linebreak
In fact, solutions are known to increase
quite substantially in size in many cases,
as shown in the examples below (Figs.\;1, 2).
In view of the constraint (1.5),
however, \linebreak
any substantial growth of
$ \,\!|\:\!u\:\!|\,\! $
leads to the development
of high frequency structures,
as illustrated below,
which in turn tend
to be efficiently dissipated
by the ever larger viscosity
present in these regions
(the local bulk viscosity
is proportional \linebreak
to
$ |\;\!u\;\!|^{\:\!\alpha} \!\;\!$).
Thus,
although
the basic ingredients for
solution blow-up
are clearly there,
especially for large $\kappa$,
the final outcome seems difficult
to predict,
be it on physical or mathematical grounds.
This interesting interaction
between convection and diffusion
due to $ \beta(x,t) > 0 $
has not been studied
in the literature
(see e.g.\;\cite{BandleBrunner1998, DengLevine2000, %
QuittnerSouplet2007}). \\
\mbox{} \vspace{-0.700cm} \\
\mbox{} \vspace{-0.350cm} \\
\mbox{} \hspace{-0.800cm}
\begin{minipage}[t]{14.600cm}
\includegraphics[keepaspectratio=true, scale=0.9]{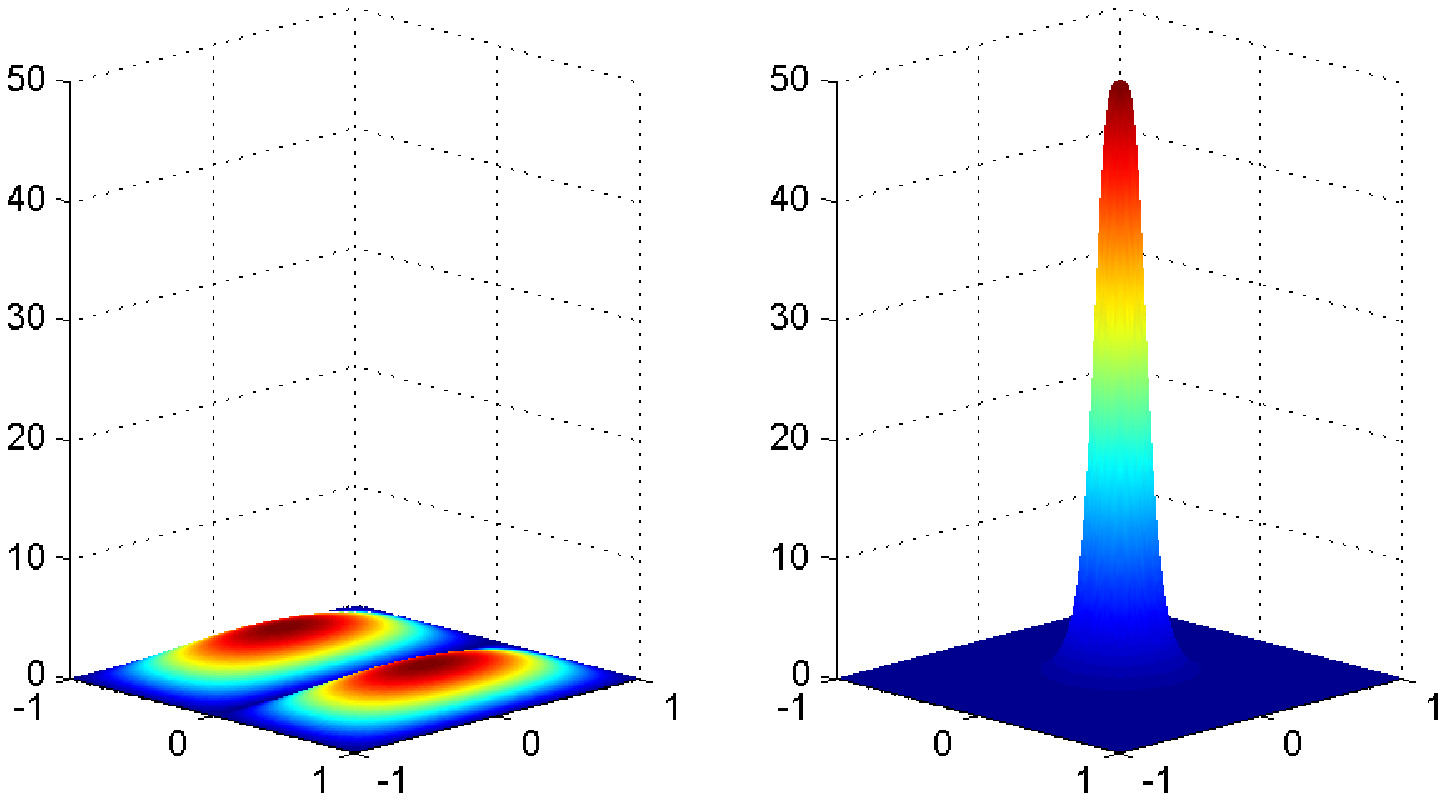}
\mbox{} \vspace{-0.550cm} \\
\mbox{} \hspace{+1.250cm}
{\small {\bf Fig.\;1:}
The solution $u(\cdot,t)$ at time $ \!\;\!\;\!t = 1000 $ (right)
for some given initial state \linebreak
\mbox{} \hspace{+1.275cm}
%
%
%
\!compactly supported in the square
%
%
$\:\!|\;\!x\;\!| \leq 1 $,
$\:\!|\;\!y\;\!| \leq 1 $
(left),
in the case $\;\! n = 2 $, \linebreak
\mbox{} \hspace{+1.275cm}
$ \:\!\alpha = \kappa = 1 \:\!$
and
$ \;\!\mbox{\boldmath $b$}(x,t) = - \,|\;\!x\;\!|^{\:\!2} \;\!x
\;\!/\!\;\!\;\!(10^{-\;\!4} +\!\;\!\;\!|\;\!x\;\!|^{\:\!4}\!\;\!/\;\!4\;\!) $,
showing an 18-fold increase \linebreak
\mbox{} \hspace{+1.275cm}
in solution size.
An impressive
18,000-fold increase
in size
would have been ob-\linebreak
\mbox{} \hspace{+1.275cm}
served
in this example
by taking
$ \;\!\mbox{\boldmath $b$}(x,t) = - \,|\;\!x\;\!|^{\:\!2} \;\!x
\;\!/\!\;\!\;\!(10^{-\;\!10} +\!\;\!\:\!|\;\!x\;\!|^{\:\!4}\!\;\!/\;\!4\;\!)
\:\!$
instead. \\
}
\end{minipage}
\mbox{} \vspace{-1.500cm} \\
\mbox{} \hspace{-0.800cm}
\begin{minipage}[t]{14.600cm}
\includegraphics[keepaspectratio=true, scale=0.9]{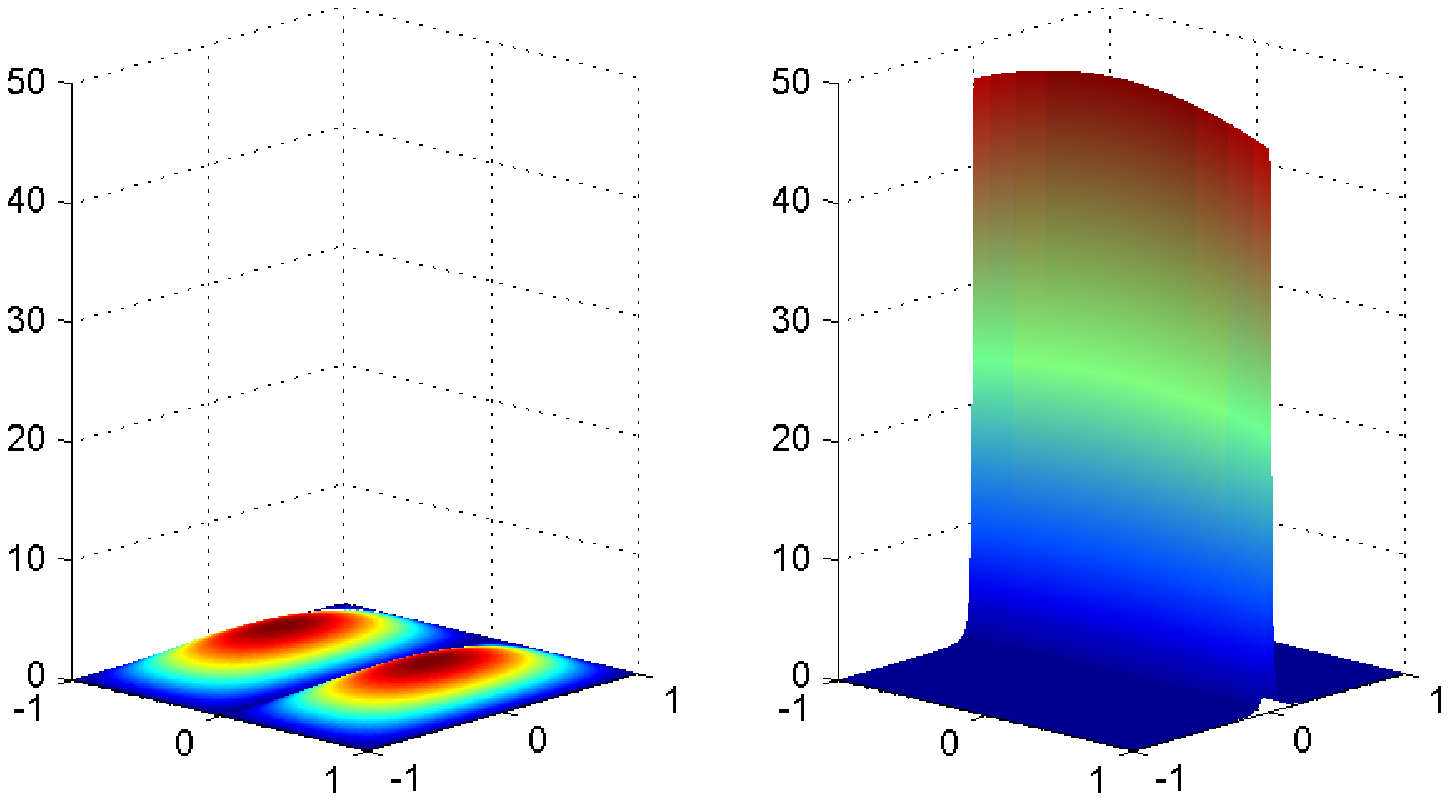}
\mbox{} \vspace{-0.650cm} \\
\mbox{} \hspace{+1.250cm}
{\small {\bf Fig.\;2:}
The solution of (1.7)
at time $ \!\;\!\;\!t = 1000 $ (right)
for the same initial state \linebreak
\mbox{} \hspace{+1.275cm}
considered in Fig.\;1,
assuming $\;\! n = 2 $,
$ \:\!\alpha = 2 $, $ \:\!\kappa = 1 \:\!$
and
$ \!\;\!\;\!\mbox{\boldmath $b$}(x,t) = (b_{1}(x), b_{2}(x)\,\!) $ \linebreak
\mbox{} \hspace{+1.275cm}
with
$ \;\!b_{1}(x) = - \;\!4/25\, x_{\mbox{}_{\scriptstyle 1}}
\!\;\!/\!\;\!\;\![\;\!
(16 + x_{1}^{\;\!2})^{\mbox{}^{\scriptstyle \!2}} \!\;\!
(10^{-\;\!4} \!\;\!+ x_{2}^{\;\!2})\!\;\!\;\!] $,
$ \!\;\!\;\!b_{2}(x) =
- \;\!4/25\, x_{\mbox{}_{\scriptstyle 2}}
/\!\;\!\;\![\;\!
(16 + x_{1}^{\;\!2}) \cdot $ \\
\mbox{} \hspace{+1.275cm}
$ (10^{-\;\!4} \!\;\!+ x_{2}^{\;\!2})^{\mbox{}^{\scriptstyle \!2}}] $,
showing once again an 18-fold increase
in solution size.
A similar \linebreak
\mbox{} \hspace{+1.275cm}
18,000-fold increase
in size
would happen
with
$ \;\!b_{1}(x) = - \;\!10^{-\;\!5}\;\! x_{\mbox{}_{\scriptstyle 1}}
\!\;\!/\!\;\!\;\![\,
(4 \!\;\!\;\!+\!\;\!\;\! x_{1}^{\;\!2})^{\mbox{}^{\scriptstyle \!2}}
\!\!\;\!\;\!\cdot $ \linebreak
\mbox{} \hspace{+1.275cm}
$ (10^{-\;\!10} \!\;\!+ x_{2}^{\;\!2})\!\;\!\;\!] \;\!$
and
$ \;\!b_{2}(x) = - \;\!10^{-\;\!5}\;\! x_{\mbox{}_{\scriptstyle 2}}
/\!\;\!\;\![\,
(4 \!\;\!\;\!+\!\;\!\;\! x_{1}^{\;\!2})
\;\!
(10^{-\;\!10} \!\;\!+ x_{2}^{\;\!2})^{\mbox{}^{\scriptstyle \!2}}]
\!\;\!\;\!$
in this example. \\
}
\end{minipage}
\nl
\mbox{} \vspace{-0.500cm} \\

Going back to the general setting
(1.1), (1.2),
we now state our main results.
Even when solutions are subject to
considerable growth
due to strong convection instabilities,
the constraint (1.5) makes
diffusion ultimately have the upper hand \linebreak
in {\em every\/} case,
thus preventing any finite time blow-up
from happening
(i.e., $\mbox{\small $T$}_{\!\:\!\ast} \!= \infty $): \\
\nl
\mbox{} \vspace{-0.775cm} \\
\mbox{} \hspace{-0.500cm}
\fbox{%
\mbox{} \hspace{-0.075cm}
\begin{minipage}[t]{14.950cm}
\nl
\mbox{} \vspace{-0.275cm} \\
%
%
%
%
{\bf Theorem I.}
\textit{%
Let
$ \;\!\kappa \geq 0$.
\!\!\;\!Then,
all solutions
to the problem $(1.1)$, $(1.2)$
satisfying $(1.3)$
are {\bf \em globally defined\/}
$($i.e, defined for all $\,t > 0\:\!)$.
} \\
\mbox{} \vspace{-0.475cm} \\
\end{minipage}
\:\!
}

\nl
\mbox{} \vspace{-0.050cm} \\
\mbox{\small \sc Theorem I} is established
after a series of technical lemmas
in
\cite{DiehlFabrisZingano2018}.
It significantly improves a previous result
obtained in \cite{Diehl2015, Fabris2013},
which was restricted to
vanishing viscosity solutions
and in addition
required the extra assumption that
$ \;\!\kappa \:\!<\:\! \alpha + 1/n $.
From the lemmata in
\cite{DiehlFabrisZingano2018}
we also obtain
an important pointwise estimate
for the solutions
of (1.1), (1.2)
involving the quantities \\
\mbox{} \vspace{-0.050cm} \\
\mbox{} \hspace{+3.500cm}
$ {\displaystyle
\mathbb{U}_{q}\!\;\!(\:\!t\:\!)
\,:=
\sup_{0 \,<\,\tau \,<\,t}
\|\, u(\cdot,\tau) \,
\|_{\mbox{}_{\scriptstyle L^{q}\!\;\!(\mathbb{R}^{n})}}
\qquad
\;\;
(\:\!1 \:\!\leq\;\! q \;\!\leq \infty\,\!)
} $
\hfill (1.8$a$) \\
and \\
\mbox{} \vspace{-1.000cm} \\
\begin{equation}
\tag{1.8$b$}
\mathbb{F}_{\!\:\!\mu}\!\;\!(\:\!t\:\!)
\,:=
\sup_{0 \,<\,\tau \,<\,t}
\frac{\;\!F(\,\!\tau)\;\!}{\mu(\,\!\tau)}
\end{equation}
\mbox{} \vspace{-0.050cm} \\
where $F$ is given in (1.2) above.
Let $\;\! a := n\:\!(\kappa - \alpha)$.
Taking $ \;\!p \geq 1 $,
$ \sigma > 1 $ satisfying \\
\mbox{} \vspace{-0.550cm} \\
\begin{equation}
\tag{1.9}
p \,>\, n\:\!(\kappa - \alpha),
\qquad
\sigma
\:\!\geq\;
\max\,\Bigl\{\,
\mbox{\small $ {\displaystyle \frac{\;\!2\,}{\mbox{\normalsize $p$}} }$},
\;
\mbox{\small $ {\displaystyle
1 \:\!+\,\frac{\:\!(2\:\!\mbox{\normalsize $\kappa$} -
\mbox{\normalsize $\alpha$})_{\mbox{}_{-}}}{\mbox{\normalsize $p$}} }$}
\,\Bigr\}
\end{equation}
\mbox{} \vspace{-0.150cm} \\
it is shown in
\cite{DiehlFabrisZingano2018}
the following fundamental estimate
(see also \cite{BarrionuevoOliveiraZingano2014, %
BrazMeloZingano2015, Diehl2015, Fabris2013, Zingano2015}): \\
\mbox{} \vspace{-0.150cm} \\
%
%
%
%
%
\mbox{} \hspace{-0.650cm}
\fbox{%
\mbox{} \;\!\;\!\;\!
\begin{minipage}[t]{14.900cm}
\mbox{} \vspace{+0.450cm} \\
{\bf Theorem II.}
\textit{%
\,Let
$ \,p \geq 1 $,
$ \sigma > 1 $
satisfy $\:\!(1.9)$.
Then \\
}
\mbox{} \vspace{-0.650cm} \\
\begin{equation}
\tag{1.10}
\mathbb{U}_{\infty}(\:\!t\:\!)
\;\!\!\;\!\;\!\leq\;\!\;\!
\mathbb{K}
\:\!\cdot\,
\max\,\biggl\{\;\!
\|\, u_0 \;\!\|_{\mbox{}_{\scriptstyle
L^{\infty}(\mathbb{R}^{n})}}
\!\;\!;
\;\!\;\!
\mathbb{F}_{\mu}\!\;\!(\:\!t\:\!)^{\mbox{}^{\scriptstyle
\!
\frac{\scriptstyle n}{\scriptstyle \;\!p \;\!-\;\! a\:\!}
}}
%
%
\mathbb{U}_{\mbox{}_{\scriptstyle \!\;\! p }}
\!\;\!(\:\!t\:\!)^{\mbox{}^{\scriptstyle
\!
\frac{\scriptstyle p }
{\scriptstyle \;\!p \;\!-\;\! a \;\!}
}}
\biggr\}
\end{equation}
\mbox{} \vspace{-0.100cm} \\
\textit{%
for every
$ \;\!t > 0 $,
and some constant
$\;\!\mathbb{K} = \mathbb{\small K}(n,\kappa,\alpha,p,\sigma) > 1 $,
\!\;\!where
$\;\! a = n\:\!(\kappa - \alpha) $. \\
\mbox{} \vspace{-0.400cm} \\
}
\end{minipage}
\;
}
%
%
\nl
\nl
\mbox{} \vspace{-0.525cm} \\
%
%

%
%

%
{\bf 2. Some generalizations} \\
\mbox{} \vspace{-0.650cm} \\

We note that the results
discussed above
apply to more general degenerate
parabolic equations of the form \\
\mbox{} \vspace{-0.600cm} \\
\begin{equation}
\tag{1.1}
u_t \,+\: \mbox{div}\;\!\mbox{\boldmath $f$}(x,t,u)
\,+\: \mbox{div}\,\mbox{\boldmath $g$}(t,u)
\:=\;
\mbox{div}\,\mbox{\boldmath $A$}(x,t,u,\nabla u).
\end{equation}
\mbox{} \vspace{-0.200cm} \\
where
$ \mbox{\boldmath $f$} \!\;\!$, $ \mbox{\boldmath $g$} \!\;\!\;\!$
are as before
and
$ {\displaystyle
\!\;\!\;\!
\mbox{\boldmath $A$} \in
C^{0}(\:\!\mathbb{R}^{n} \!\times\!\;\! [\;\!0, \:\!\infty)
\!\!\;\!\;\!\times \mathbb{R}
\!\;\!\times\!\;\!\mathbb{R}^{n})
} $
satisfies
the condition \\
\mbox{} \vspace{-0.525cm} \\
\begin{equation}
\tag{2.2}
\langle \;\!\mbox{\boldmath $A$}(x,t,\mbox{u},{\bf v}),
\;\!{\bf v} \;\!\rangle
\;\!\geq\;\!\;\!
\mu(t) \,|\;\!\mbox{u}\;\!|^{\:\!\alpha}\;\!
|\,{\bf v}\,|^{\:\!2}
\quad
\mbox{and}
\quad
\;
|\, \mbox{\boldmath $A$}(x,t,\mbox{u},{\bf v}) \;\!|
\;\!\;\!\leq\;\!\;\!
M(t) \;\!\;\!|\;\!\mbox{u}\;\!|^{\:\!\alpha}
\;\!
|\,{\bf v}\,|
\end{equation}
\mbox{} \vspace{-0.150cm} \\
for all
$ {\displaystyle
\:\!
x \in \mathbb{R}^{n}\!\:\!,
\;\!\;\!
t \:\!\geq\:\!0, \,
\mbox{u} \in \mathbb{R}, \,
{\bf v} \in \mathbb{R}^{n}
} $\!\:\!,
\!\;\!\;\!and some
$ {\displaystyle
\;\!
\mu, \!\;\!\;\!M
\in C^{0}(\,\![\;\!0, \:\!\infty))
} $,
\:\!with
$ \:\!\mu(t) > 0 $, \linebreak
as long as
it can be shown
that
$ {\displaystyle
\;\!
\|\, u(\cdot,t) \,
\|_{{\scriptstyle L^{1}(\mathbb{R}^{n})}}
\!\;\!
} $
cannot blow up in finite time.
The results also extend
to the case of arbitrary initial states
$ \;\!u_0 \!\;\!\in L^{1}(\mathbb{R}^{n}) $
\mbox{[\,}not necessarily bounded\;\!\mbox{]},
with condition (1.3$b$)
then replaced by the assumption
that
$ \;\!u(\cdot,t) \in L^{\infty}_{\mbox{\scriptsize loc}}
(\,\!(\:\!0, \:\!\mbox{\small $T$}_{\!\ast}), \:\!
L^{\infty}(\mathbb{R}^{n})\,\!) $,
or even more generally
to initial data
$ \;\!u_0 \!\;\!\in L^{p_{\mbox{}_{0}}}\!\!\;\!\;\!(\mathbb{R}^{n}) $
for some given $\;\!1 \!\;\!\leq p_{\mbox{}_{0}} \!< \infty $,
\:\!with only minor changes in the statements,
provided once more that it can be shown that
$ {\displaystyle
\|\, u(\cdot,t) \,
\|_{{\scriptstyle L^{q}(\mathbb{R}^{n})}}
\!\;\!
} $
will not blow up in finite time
for some suitable
$ \;\!p_{\mbox{}_{0}} \!\;\!\leq q < \infty $.
\newpage
%
%
%
\mbox{} \vspace{-0.950cm} \\

{\bf Acknowledgements.}
This work was partially supported
by {\small CAPES}
(\mbox{\small M}inistry of \mbox{\small E}ducation,
\mbox{\small B}razil),
\mbox{\small G}rant
\mbox{\footnotesize \#}\,\mbox{\small 88887.125079/2015}.
The computations in this research
were performed by the {\small \sc SGI} cluster
{\small \sc Altix} 1350/450 of the
{\small \sc Centro Nacional de Processamento de Alto
Desempenho em S\~ao Paulo (CENAPAD-SP)},
Brazil. \\

%
%

%
\mbox{} \vspace{-1.350cm} \\
%
%
%
%
%
%

%
%

%
\nl
\mbox{} \vspace{-0.050cm} \\
\nl
{\small
\begin{minipage}[t]{10.00cm}
\mbox{\normalsize \textsc{Nicolau Matiel Lunardi Diehl}} \\
Instituto Federal de Educa\c c\~ao, Ci\^encia e Tecnologia \\
Canoas, RS 92412, Brazil \\
E-mail: {\sf nicolau.diehl@canoas.ifrs.edu.br} \\
\end{minipage}
\nl
\mbox{} \vspace{-0.450cm} \\
\nl
\begin{minipage}[t]{10.00cm}
\mbox{\normalsize \textsc{Lucin\'eia Fabris}} \\
Coordenadoria Acad\^emica \\
Universidade Federal de Santa Maria \\
Campus de Cachoeira do Sul \\
Cachoeira do Sul, RS 96501, Brazil \\
E-mail: {\sf lucineia.fabris@ufsm.br} \\
\end{minipage}
\nl
\mbox{} \vspace{-0.450cm} \\
\nl
\begin{minipage}[t]{10.00cm}
\mbox{\normalsize \textsc{Paulo Ricardo de Avila Zingano}} \\
Departamento de Matem\'atica Pura e Aplicada \\
Universidade Federal do Rio Grande do Sul \\
Porto Alegre, RS 91509, Brazil \\
E-mail: {\sf paulo.zingano@ufrgs.br} \\
\mbox{} \hspace{+1.150cm}
        {\sf zingano@gmail.com} \\
\end{minipage}
}
%
%

\end{document}